\documentclass[12pt]{article}
\usepackage[]{latexsym}

\newcommand{\rit}{{\rm I\!R}}
 
\newcommand{\pre}{{\bf Proof.\ }}

\newcommand{\cs}{{\mathcal S}}
\newcommand{\ce}{{\mathcal E}}

\newcommand{\cd}{{\mathcal D}}
\newcommand{\ct}{{\mathcal T}}

\newcommand{\cf}{{\mathcal F}}

\newcommand{\tr}{{\rm tr}} 

\newtheorem{theo}{Theorem}[section]
\newtheorem{prop}[theo]{Proposition}
\newtheorem{lema}[theo]{Lemma}

\title{On the cohomology of the spaces of differential operators acting on skewsymmetric tensor fields or on forms, as modules of the Lie algebra of vector fields}
\author { B. Agrebaoui, F. Ammar\\
{\small (Faouzi.Ammar@mail.rnu.tn, B.Agreba@fss.rnu.tn) }\\
{\small   Department of Mathematic, Faculty of Sciences of Sfax,}\\
{\small Route Soukra  BP 802, 3018 Sfax Tunisia} \\
P. Lecomte\\
{\small (plecomte@ulg.ac.be)}\\
{\small University of Liege, Department of Mathematics,}\\
{\small Grande Traverse, 12, B37, B4000 Sart Tilman}}
\date{}
 
\begin{document}
\maketitle
\abstract{The spaces of differential operators acting on skewsymmetric tensor fields or on smooth forms of a smooth manifold are representations of its Lie algebra of vector fields. We compute the first cohomology spaces of these representations and show how they are related to the cohomology with coefficients in the space of smooth functions of the manifold.}
\section{Introduction} 
In this paper, we will compute the first cohomology spaces of the Lie algebra of vector fields of a smooth manifold with coefficients  in the space of linear differential operators acting on contravariant skewsymmetric tensor fields or on differential forms of the manifold.

As we shall explain in more details below, we show that these spaces are generated by cohomology classes canonically constructed out of the cohomology of the Lie algebra of vector fields acting on the space of smooth functions (the latter is computed in \cite{8}).

Our computation is the first step towards the study of the deformations of the considered modules. Using the results of the present paper, the miniversal deformation has been computed in \cite{2} in the case of forms. It is shown that it is equivalent to the infinitesimal one, a less rich situation than that of the symmetric case.

We will follow mostly the same strategy than in \cite{6}, where the operators act on symmetric contravariant tensor fields instead on skewsymmetric tensor fields. The method goes as follows. One studies the locality of $0$- and $1$-cocycles. One then performs the computation when $M=\rit^m$ and then extend it to arbitrary $M$ by some standard gluing arguments. When $M=\rit^m$, we filter the cohomology with the projective embbeding $sl_{m+1}$ of $sl(m+1,\rit)$ as a Lie subalgebra of $Vect(\rit^m)$. This allows to use the results of \cite{5}. We then extend to the whole algebra of vector fields using relative cohomology.

As far as local computations are concerned, the calculations are essentially the same for the two representations we are dealing with. Therefore, we detail explicitely only the case of contravariant skewsymmetric tensor fields, indicating just how the results adapt for the differential forms.

The computations are valid for manifolds of dimension at least $2$.

\section{Definitions - Notations}\label{Def}

\paragraph{a)}Let $M$ be a connected, second countable Hausdorff manifold of class $C^\infty$.
Denote by $Vect(M)$, $\wedge(M)$ and $\Omega(M)$, the Lie algebra of smooth vector fields, the space of smooth skewsymmetric tensor fields and that of smooth forms of $M$ respectively. The space $\wedge(M)$ is graded by the tensor degree: $\wedge(M)=\bigoplus_i\wedge^i(M)$. Each homogeneous component is a representation of $Vect(M)$ that acts on it by Lie derivative. Similarly, we consider $\Omega(M)=\bigoplus_i\Omega_i(M)$ as a graded module of $Vect(M)$, the gradings being induced by the degree of forms.

For any two spaces of sections $\ce$ and $\cf$ of vector bundles over $M$, we denote by $\cd^k(\ce,\cf)$ the space of $k$th-order linear differential operators from $\ce$ into $\cf$. This gives a filtration of the space $\cd(\ce,\cf)=\bigcup_k\cd^k(\ce,\cf)$ of all linear differential operators from $\ce$ into $\cf$. If $\ce=\cf$ then we will often write simply $\cd^k(\ce)$ and $\cd(\ce)$ instead of $\cd^k(\ce,\cf)$ and $\cd(\ce,\cf)$ respectively. 

If $Vect(M)$ acts on $\ce$ and $\cf$, then it also acts on $\cd(\ce,\cf)$ in the standard way 
\[
L_XA=L_X\circ A-A\circ L_X.
\]
In most of the cases, this actions also preserves the filtrations. It is the case when $\ce$ and $\cf$ are one of the representations $\wedge (M)$ or $\Omega(M)$.

Our goal is to compute the spaces $H^i(Vect(M),\cd^k(\wedge^p(M),\wedge^q(M)))$ and $H^i(Vect(M),\cd^k(\Omega_p(M),\Omega_q(M)))$, $i=0,1$.

\paragraph{b)} Recall that, on $\rit^m$, any linear $k$th-order differential operator can be writen under the form
\[
\sum_{r\leq k}\sum_{i_1\cdots i_r}A^{i_1\cdots i_r}\circ \partial_{i_1\cdots i_r},
\]
where $A^{i_1\cdots i_r}$ belongs to $C^\infty(\rit^m,V^*\otimes W)$,  $V$ and $W$ being the typical fibers of the bundles of which the arguments and the values of the operator are sections.
Moreover
\[
\partial_{i_1\ldots i_r}=\frac{\partial}{\partial x^{i_1}}\cdots\frac{\partial}{\partial x^{i_r}}
\]
acts on $C^\infty(\rit^m,V)$ by ordinary partial derivatives.

In typical situations, $V$ and $W$ are representations of $gl(m,\rit)$ and the corresponding Lie derivative of $Vect(\rit^m)$ on $C^\infty(\rit^m,V)$ is just given by
\[
L^\rho_Xf=X.f-\rho(DX)f.
\]
(Here, $X.f$ is the usual derivative of $f$ in the direction of $X$, $DX$ is the Jacobian matrix of $X$ and $\rho$ denotes the action of $gl(m,\rit)$ on $V$.)

\paragraph{c)}
The Lie algebra $sl(m+1, \rit)$ of matrices of vanishing trace, has a decomposition of the shape :
\[
\rit^m \bigoplus gl(m, \rit)\bigoplus {\rit^m}^*
\]
where for $h \in \rit^m$ and $\alpha \in {\rit^m}^*$, the braket is given by
$[h,\alpha]=\alpha(h)1+h\bigotimes \alpha$, the other brackets are the obvious one.
This decomposition gives $sl(m+1, \rit)$ a structure of graded Lie algebra where $\rit^m$, 
$gl(m, \rit)$ and ${\rit^m}^*$ are the homogeneous components of degree $-1, \ 0$
and $1$ respectively. We can realize $sl(m+1, \rit)$ as a subalgebra of  $Vect(\rit^n)$ where
the vector fields associated to $h=(h^i)\in \rit^m$, $A=(A^i_j)\in gl(m, \rit)$ and  
$\alpha=(\alpha_i)\in {\rit^m}^*$ are respectively
\[
h^*=-h^i\partial_i, \ A^*=-A^i_jx^j\partial_i \mbox{ and } \alpha^*=\alpha(x)x^i\partial_i
\]
For the sake of brevity, we denote $sl_{m+1}$ the above realization of $sl(m+1, \rit)$.

We will need the following result, taken from \cite{5}.
Let the map 
\[
\chi: \wedge(gl(m,\rit),\wedge(\rit^{m*},V))\to\wedge(sl_{m+1},C_\infty(\rit^m,V))
\]
be given by
\[
\begin{array}{l}
(X_0,\ldots,X_{t+u-1})\to\\[1ex]
\frac{(-1)^t}{t!u!(m+1)^u}\sum_\nu sign(\nu)(\gamma(DX_{\nu_0},\ldots,DX_{\nu_{t-1}})(d\tr (DX_{\nu_t}),\ldots,d\tr (DX_{\nu_{t+u-1}}))). 
\end{array}
\]
Then
\begin{theo}\label{SlmCohom}
The map induced in cohomology
\[
\chi_\sharp: H(gl(m,\rit),\wedge(\rit^{m*},V))\to H(sl_{m+1},C_\infty(\rit^m,V))
\]
is a bijection.
\end{theo}

\section{The space $H^k(sl_{m+1},\cd(\wedge^p(\rit^m),\wedge^q(\rit^m)))$}
To compute that cohomology space, we proceed by induction on $k$, like in \cite{5}, using the short exact
sequence 
\[
 0\rightarrow \cd^{k-1}(\wedge^p,\wedge^q)\stackrel{i}{\rightarrow} \cd^k(\wedge^p,\wedge^q)\stackrel{\sigma}{\rightarrow} \cs^k (\wedge^p,\wedge^q) \rightarrow 0
\]
that induces
an exact triangle
\begin{equation} \label{TriangleExact}
\begin{array}{cl}
 H(sl_{m+1},\cd^{k-1}(\wedge^p,\wedge^q))\\
& \searrow i_\sharp \\
\theta\;\big\uparrow & H(sl_{m+1},\cd^k(\wedge^p,\wedge^q))  \\
& \swarrow \sigma_\sharp \\
H(sl_{m+1},\cs^k(\wedge^p,\wedge^q))
\end{array}
\end{equation}
To simplify the notations, we have replaced $\wedge^*(\rit^m)$ by $\wedge^*$. Moreover, we denote by $\cs^k(\wedge^p,\wedge^q)$ the space
\[
C^\infty(\rit^m,\vee^k\rit^m\otimes\wedge^p\rit^{m*}\otimes\wedge^q\rit^m)
\]
of $k$-symmetric contravariant tensor fields valued in $Hom(\wedge^p,\wedge^q)$, which is isomorphic to the space of principal symbols of $\cd^k(\wedge^p,\wedge^q)$. The map $\theta$ is the induced connecting homomorphism; it is of degree $1$ \cite{3}.

\begin{prop}\label{CohomSymb}
a)If $k\geq 2$ or if $p<q$ then $H(sl_{m+1}, \cs^k(\wedge^p,\wedge^q))=0$.\\[1ex]
b)The space $H(sl_{m+1}, \cs^0(\wedge^p,\wedge^q))$ is isomorphic to  $(\wedge gl(m,\rit)^*)_{g-inv}$.\\[1ex]
c)The space $H(sl_{m+1}, \cs^1(\wedge^p,\wedge^q))$ is isomorphic to  $(\wedge gl(m,\rit)^*)_{g-inv}$ if $p>q$ and is vanishing if $p=q$.
\end{prop}
\pre It follows from Theorem \ref{SlmCohom} that, if it is not vanishing, the space $H^u(sl_{m+1}, \cs^k(\wedge^p,\wedge^q))$ is isomorphic to
\[
H^{u-p+q+k}(gl(m,\rit),\wedge^{p-q-k}(\rit^{m*},V)),
\]
where $V=\vee^k\rit^m\otimes\wedge^p\rit^{m*}\otimes\wedge^q\rit^m$ is equipped with the canonical action of $gl(m,\rit)$. By \cite{5}, Proposition 4.3, 
\[
\begin{array}{l}
H^{u-p+q+k}(gl(m,\rit),\wedge^{p-q-k}(\rit^{m*},V))=\\[1ex]
(\wedge^{u-p+q+k}gl(m,\rit)^*)_{g-inv}\bigotimes (\wedge^{p-q-k}(\rit^{m*},V))_{s-inv}
\end{array}
\]
where $g-inv$ and $s-inv$ denotes the invariant elements with respect to $gl(m,\rit)$ and $sl(m,\rit)$ respectively. It follows then easily from the theory of invariants of classical groups (see e.g. \cite{4}) that $(\wedge^{p-q-k}(\rit^{m*},V))_{s-inv}$ is non vanishing only if $k=0$ or $k=1$. In both cases, it is 1-dimensional. It is spanned by the mapping
\[
I_0:(\alpha_1,\ldots,\alpha_{p-q})\mapsto (T\mapsto i_{\alpha_1}\cdots i_{\alpha_{p-q}}T)
\]
in the first case and by
\[
I_1:(\alpha_1,\ldots,\alpha_{p-q-1})\mapsto (T\mapsto i_\eta i_{\alpha_1}\cdots i_{\alpha_{p-q-1}}T)
\]
in the second case. (Here, we view an element of $V$ as being an homogeneous polynomial of degree $k$ in $\eta\in \rit^{m*}$ valued in the space of linear mappings from $\wedge^p \rit^m$ into $\wedge^q \rit^m$.) The result then immediately follows.

\begin{lema}\label{Connecting}
Assume that $k=1$ and $p>q$ and identify the source and target of the connecting homomorphism in (\ref{TriangleExact}) with $(\wedge gl(m,\rit)^*)_{g-inv}$. One has
\[
\theta(\gamma)=(-1)^{|\gamma|}(p-q+1)(m+1)\gamma, \forall\gamma\in (\wedge gl(m,\rit)^*)_{g-inv}.
\]
\end{lema}
\pre (sketch) Let $\gamma\in (\wedge^a gl(m,\rit)^*)_{g-inv}$ be given. According to the above proof, the coboundary $\partial\chi(\gamma\otimes I_1)$ of $\chi(\gamma\otimes I_1)$ is of the form $ \chi(\gamma'\otimes I_0)$ and we have to show that $\gamma'=(-1)^a(p-q+1)(m+1)\gamma$. Set $b=p-q$. Since
\[
\chi(\gamma'\otimes I_0)(A^*_0,\ldots,A^*_{a-1},\alpha^*_0,\ldots,\alpha^*_b)=\gamma'(A_0,\ldots,A_{a-1})I_0(\alpha_0,\ldots,\alpha_b)
\]
where $\gamma'$ and $I_0$ have constant coefficient, it suffices to compute
\begin{equation}\label{cobord}
\partial\chi(\gamma\otimes I_1)(A^*_0,\ldots,A^*_{a-1},\alpha^*_0,\ldots,\alpha^*_b)_{|_{x=0}}
\end{equation}
Recall that the coboundary $(\partial c)(X_0,\ldots,X_t)$ writes
\[
\sum_i (-1)^iL_{X_i}c(X_0,\ldots,\hat{X_i},\ldots,X_t)+\sum_{i<j}(-1)^{i+j}c([X_i,X_j],\ldots,\hat{X_i},\ldots,\hat{X_j},\ldots).
\]
One easily sees that, in (\ref{cobord}), the terms corresponding to the second sum are vanishing as well as, in the first sum, these for which $X_i$ is one of the $A^*_j$. One is thus left to evaluate
\[
\begin{array}{l}
\sum_i (-1)^{a+i}L_{\alpha_i^*}\chi(\gamma\otimes I_1)(A^*_0,\ldots,A^*_{a-1},\alpha^*_0,\ldots,\hat{\alpha_i^*},\ldots,\alpha^*_b)_{|_{x=0}}=\\[1ex]
(-1)^a\gamma(A_0,\ldots,A_{a-1})\sum_i (-1)^iL_{\alpha_i^*}I_1(\alpha^*_0,\ldots,\hat{\alpha_i^*},\ldots,\alpha^*_b).
\end{array}
\]
A direct computation shows that 
\[
(-1)^iL_{\alpha_i^*}I_1(\alpha^*_0,\ldots,\hat{\alpha_i^*},\ldots,\alpha^*_b)_{|_{x=0}}=(m+1)I_0(\alpha_0,\ldots,\alpha_b).
\]
Hence the Lemma.

The next theorem imediately follows from Proposition \ref{CohomSymb} and Lemma \ref{Connecting}

\begin{theo} One has\\
(a) If $p<q$, then $H(sl_{m+1}, \cd^k(\wedge^p,\wedge^q))=0$ for all $k$.\\[1ex]
(b) If $p=q$, then the spaces $H(sl_{m+1}, \cd^k(\wedge^p,\wedge^q))$, $k\geq 0$, and $H(sl_{m+1}, \cd(\wedge^p,\wedge^q))$ are isomorphic to $\wedge (gl(m,\rit)^*)_{g-inv}$.\\[1ex]
(c) If $p>q$, then $H(sl_{m+1}, \cd^0(\wedge^p,\wedge^q))$ is isomorphic to $\wedge (gl(m,\rit)^*)_{g-inv}$; the spaces $H(sl_{m+1}, \cd^k(\wedge^p,\wedge^q))$, $k\geq 1$, and $H(sl_{m+1}, \cd(\wedge^p,\wedge^q))$ are vanishing.\\
In particular
\[
H^1(sl_{m+1}, \cd(\wedge^p,\wedge^q))=
\left\lbrace
\begin{array}{l}
\rit \mbox { if } p=q\\
0 \mbox { otherwise }
\end{array}\right.
\]
\end{theo}

\section{The space $H^0(Vect(M),\cd^k(\wedge^p(M),\wedge^q(M)))$}
We use the notations of Section \ref{Def}. The dimension of the manifold $M$ is at least $2$.

\begin{theo}\label{0Cohom}For every $k\geq 0$, one has
\[
H^0(Vect(M), \cd^k(\wedge^p(M),\wedge^q(M)))=
\left\{
\begin{array}{l}
\rit \mbox { if } p=q\\
0 \mbox { otherwise }
\end{array}
\right.
\]
\end{theo}
\pre Let $\ct\in\cd^k(\wedge^p(M),\wedge^q(M)$ be $Vect(M)$-equivariant. It is a local map. Indeed, assume that $T\in\wedge^p(M)$ vanishes on some open subset of $M$ and let $a$ be an arbitray point of that subset. According to \cite{10}, there exist finitely many vector fields $X_i$ and tensors $T_i$ vanishing in a neighborhood of $a$ such that
\[
T=\sum L_{X_i}T_i.
\]
Therefore
\[
\ct(T)=\sum L_{X_i}(\ct(T_i))
\]
vanishes in that neighborhood. Hence, $\ct$ is local.
It follows from a well known theorem of Peetre that the restriction of $\ct$ over any relatively compact domain of chart $U$ of $M$ is a differential operator. It is moreover $sl_{m+1}$-equivariant. Applying theorem \ref{SlmCohom}, we see that $\ct_{|U}$ vanishes if $p\neq q$ or is some constant multiple $k_U$ of the identity otherwise. Clearly, $M$ being conneced, $k_U$ is independant of $U$. The result then follows immediately.

\section{The space $H^1(Vect(M),\cd^k(\wedge^p(M),\wedge^q(M)))$}
We use again the notations of Section \ref{Def}. The dimension of the manifold $M$ is at least $2$.

\begin{theo}\label{Cohom}One has\\
(a)  For each $q\geq 0$, $H^1(Vect(M),\cd^0(\wedge^{q+1}(M),\wedge^q(M)))\equiv\rit$.\\[1ex]
(b) For each $p,k\geq 0$, $H^1(Vect(M),\cd^k(\wedge^p(M),\wedge^p(M)))\equiv\rit\oplus H_{DR}^1(M)$.\\[1ex]
(c) In the other cases, $H^1(Vect(M),\cd^k(\wedge^p(M),\wedge^q(M)))=~0$.
\end{theo}
\pre  (sketch) Let $c:Vect(M)\to\cd^k(\wedge^p(M),\wedge^q(M))$ be a $1$-cocycle. A straightforward adaptation of the argument of \cite{6}, p. 98, shows that it is a local map. As above for $0$-cocycle, we first determine its restriction $c_{|U}$ over a relatively compact domain of chart $U$ of $M$.

It is of course a $1$-cocycle of the embedding $sl_{m+1}$ associated to the chart. From Theorem \ref{SlmCohom}, the restriction of $c_{|U}$ to $sl_{m+1}$ is of the form
\begin{equation}\label{LocCocycle}
\left\{
\begin{array}{ll}
X\to \partial b_U &\mbox{ if } p<q \mbox{ or if } p>q \mbox{ and } k>0,\\
X\to r_U \tr DX+\partial b_U(X) &\mbox{ if } p=q,\\
X\to r_U \iota_{d\tr DX}+ \partial b_U(X) &\mbox{ if } p>q  \mbox{ and } k=0,
\end{array}
\right.
\end{equation}
for some $b_U\in\cd^k(\wedge^p(U),\wedge^q(U))$ and some constant $r_U$. In each case, that mapping extends to the algebra $Vect(U)$ as a $1$-cocycle of that algebra in an obvious way.
Substrating it from $c_{|U}$, we are left with a $1$-cocycle of $Vect(U)$ that vanishes on $sl_{m+1}$.

As a mapping from $Vect(U)\times\wedge^p(U)$ into $\wedge^q(U)$, such a cocycle is a $sl_{m+1}$-equivariant differential operator. It follows first that it has constant coefficients and that it is $gl(m,\rit)$-equivariant. As easily seen, this, together with the fact that it vanishes on $sl_{m+1}$, implies that it is in fact equal to $0$. This means that $c_{|U}$ is of the form (\ref{LocCocycle}) over the whole $Vect(U)$. In other words, Theorem \ref{Cohom} is proven when $M=\rit^m$. A simple adaptation of the argument of \cite{6} allows to extend it to $M$. Hence the theorem.

\medskip
Let us describe some generators of the  cohomology space
\[
 H^1(Vect(M),\cd^k(\wedge^p(M),\wedge^q(M))).
\]
There is a canocical map
\[
H^1(Vect(M),C^\infty(M))\to  H^1(Vect(M),\cd^k(\wedge^p(M),\wedge^p(M))).
\]
It is induced by the map $c\mapsto c.id$ at the level of cocycles. From $(b)$ of Theorem \ref{Cohom}, we see that it is an isomorphism.
There is also a map
\[
H^1(Vect(M),C^\infty(M))\to  H^1(Vect(M),\cd^k(\wedge^{q+1}(M),\wedge^q(M))),
\]
induced by the correspondance that maps  a $C^\infty(M)$-valued cocycle $c$ onto the cocycle
\begin{equation}\label{Cocycle}
X\in Vect(M)\mapsto \iota(dc(X))\in\cd^k(\wedge^{q+1}(M),\wedge^q(M))).
\end{equation}
Due to the formula
\[
\iota_{d\iota_X\xi}T=L_X(\iota_\xi T)- \iota_\xi L_XT,
\]
this cocycle is a coboudary when $c$ is a closed $1$-form $\xi$, so that the map induced in cohomology could only be non trivial on the component $\rit$ of $H^1(Vect(M),C^\infty(M))$, that is spanned by any  divergence operator. It follows from the above theorem that it is only non trivial for $k=0$. This means that the cocycle (\ref{Cocycle}) is the coboudary of a differential operator of order at least $1$. Let us show directly that it is indeed the case. Let $c$ be any  $C^\infty(M)$-valued cocycle of $Vect(M)$ of order $1$. There exists an atlas of $M$ in each chart $(U,\varphi)$ of which
\[
c(X)=r\tr(DX), \forall X\in Vect(M),
\]
for some nonzero real number $r$, independant of the chart. As easily seen, $X\mapsto r\tr(DX)$ is the coboudary of the map
\[
T\mapsto r\sum_i\iota_{dx^i}\partial_iT.
\]

Due to Theorem \ref{0Cohom}, this map doesn't depend on the chart $(U,\varphi)$: it is the restriction on $U$ of a globally defined map of which $c$ is the coboundary.

\section{Cohomology with coefficients in $\cd^k(\Omega_p(M),\Omega_q(M))$}
In this section, we describe the two first cohomology spaces of $Vect(M)$ acting by Lie derivatives on the space of differential operators from the space $\Omega_p(M)$ of smooth $p$-forms into $\Omega_q(M)$. We will not give the details of the computation, just pointing out a few remarks. It is indeed exactly taylored on the same scheme than the above calculation and most of it is a straight adaptation of what we have done up to here.

We first assume that $M=\rit^m$ and filter by the subalgebra $sl_{m+1}$. We use again the short exact sequence associated to the order $k$ of differentiation to get an exact triangle analoguous to (\ref{TriangleExact}). At this stage, it is worth noticing that
\[
H(sl_{m+1},\cs^k(\wedge^p,\wedge^q))\equiv H(sl_{m+1},\cs^k(\Omega_q,\Omega_p)).
\]
Indeed, both are computed via Theorem \ref{SlmCohom} whith the representation $V$ equal to $V_{pq}=\vee^k\rit^m\otimes\wedge^p\rit^{m*}\otimes\wedge^q\rit^m$ for the lhs and $V_{qp}$ for the rhs.
In particular, the rhs is vanishing for $k>1$. Moreover, the space $(\wedge^{q-p-k}(\rit^{m*},V_{qp}))_{s-inv}$ is generated by
\[
J_0:(\alpha_1,\ldots,\alpha_{q-p})\mapsto (\omega\mapsto \alpha_1\wedge\ldots\wedge\alpha_{q-p}\wedge\omega)
\]
for $k=0$ and by
\[
J_1:(\alpha_1,\ldots,\alpha_{q-p})\mapsto (\omega\mapsto \eta\wedge\alpha_1\wedge\ldots\wedge\alpha_{q-p}\wedge\omega)
\]
for $k=1$.
In the latter case, one still has to compute the connecting homomorphism $\theta$. This is more easy than in the contravariant case and one finds that it vanishes. This leads to the following result.

\begin{theo} One has, replacing $\Omega^*(\rit^m)$ by $\Omega^*$ for simplicity,\\
(a) If $p>q$, then $H(sl_{m+1}, \cd^k(\Omega^p,\Omega^q))=0$ for all $k$.\\[1ex]
(b) If $p=q$, then the spaces $H(sl_{m+1}, \cd^k(\Omega^p,\Omega^q))$, $k\geq 0$, and $H(sl_{m+1}, \cd(\Omega^p,\Omega^q))$ are isomorphic to $\wedge (gl(m,\rit)^*)_{g-inv}$.\\[1ex]
(c) If $p<q$, then $H(sl_{m+1}, \cd^0(\Omega^p,\Omega^q))$ is isomorphic to $\wedge (gl(m,\rit)^*)_{g-inv}$; the spaces $H(sl_{m+1}, \cd^k(\Omega^p,\Omega^q))$, $k\geq 1$, and $H(sl_{m+1}, \cd(\Omega^p,\Omega^q))$ are isomorphic to $\wedge (gl(m,\rit)^*)_{g-inv}\oplus\wedge (gl(m,\rit)^*)_{g-inv}$ .\\
In particular
\[
H^1(sl_{m+1}, \cd(\Omega^p,\Omega^q))=
\left\lbrace
\begin{array}{l}
\rit \mbox { if } p=q \mbox { or } q=p+2\\
\rit^2 { if } q=p+1\\
0 \mbox { otherwise }
\end{array}\right.
\]
\end{theo}

The space $H^1(sl_{m+1}, \cd(\Omega^p,\Omega^q))$ is generated by the classes of the cocycles defined by $c_0(X):\omega\mapsto \tr(DX) \omega$, $c_{01}(X):\omega\mapsto\tr(DX) d\omega$, $c_{10}(X):\omega\mapsto d\tr(DX)\wedge\omega$ and $c_2(X):\omega\mapsto d\tr(DX)\wedge d\omega$.

These cocycles extend to the algebra $Vect(\rit^m)$ and, as in the contravariant case, one sees that they generate the space $H^1(Vect(\rit^m), \cd(\Omega^p,\Omega^q))$.

One goes from $\rit^m$ to the arbitrary manifold $M$ (still of dimension at least $2$) in the same way than in the contravariant case either. The result is thus stated below whithout proof. It can be again interpreted using the cohomology of $Vect(M)$ acting on $C^\infty(M)$, as follows: a $1$-cocycle $\gamma: Vect(M)\to C^\infty(M)$ leads to four cocycles of our cohomology, namely these given by
\[
\begin{array}{lcl}
q=p&:&\gamma(X)\omega\\
q=p+1&:&\gamma(X)d\omega\\
q=p+1&:&d(\gamma(X))\wedge\omega\\
q=p+2&:&d(\gamma(X))\wedge d\omega
\end{array}
\]

One easily sees that if a closed $1$-form is added to $\gamma$, then the cohomology class of the last two cocycles is not modified. This is not the case for the two first. More precisely, we can state
\begin{theo}One has\\
(a) For each  $q,k\geq 0$, $H^1(Vect(M), \cd^k(\Omega_q(M),\Omega_q(M)))\equiv \rit\oplus H^1_{DR}(M)$.\\[1ex]
(b) For each $q>0$, $ H^1(Vect(M), \cd^0(\Omega_q(M),\Omega_{q-1}(M)))\equiv \rit$.\\[1ex]
(c) For each $q,k>0$, $ H^1(Vect(M), \cd^k(\Omega_q(M),\Omega_{q-1}(M)))\equiv \rit^2\oplus H^1_{DR}(M)$.\\[1ex]
(d) For each $q>1, k>0$, $H^1(Vect(M), \cd^k(\Omega_q(M),\Omega_{q-2}(M)))\equiv \rit$.\\[1ex]
(e) In the other cases, $H^1(Vect(M), \cd^k(\Omega_p(M),\Omega_q(M)))$ is vanishing.
\end{theo}

Applying the previous results, we could also compute the  $0$-cohomology space of $Vect(M)$ acting on $\cd(\Omega_p,\Omega_q)$, at least that of the complex of local cochains. It's trivial to see that a $Vect(M)$-equivariant operators from $\Omega(M)$ into $\Omega_q(M)$ is local when $q>0$. This is not true for $q=0$ (one can integrate forms of maximum degree, at least if $M$ is oriented).
This would give back well known results \cite{9}, so that we will not achieve these additional computations.

\end{document}